\documentclass{article}
\usepackage{amsmath,amsthm,amsfonts}
\usepackage{amssymb}
\usepackage{epsfig}
\def\rank{\mathop{rank}}

\def\e{\epsilon}
\def\al{\alpha}
\def\la{\lambda}
\def\A{{\cal A}}
\def\B{{\cal B}}

\def\l{{l}}
\def\Imp{\mathop{Im}}

\def\C{{\bf C}}

\def\S{{\cal S}}

\def\R{{\bf R}}
\def\Z{{\bf Z}}
\def\RS{{\cal R}}
\def\spectr{\mathop{sp}}

\def\diag{\mathop{diag}}
\def\dc{\mathop{dc}}
\def\dim{\mathop{dim}}

\def\set{\mathop{set}}
\newcommand{\N}{{\bf N}}
\usepackage[latin1]{inputenc}

\newtheorem{lemma}{Lemma}

\newtheorem{theorem}{Theorem}
\newtheorem{proposition}{Proposition}
\newtheorem{corollP}{Corollary}
\theoremstyle{remark}
\newtheorem{remark}{Remark}
\title{On low rank perturbations of complex matrices and some discrete metric spaces.}
\author{Lev Glebsky and Luis Manuel Rivera}

\begin{document}
\maketitle
\section{Introduction}
The article is devoted to different aspects of the question: "What
can be done with a complex-valued matrix by a low rank
perturbation?"\footnote{The authors would like to thank
V.S.Savchenko, J. Moro and F. Dopico 
for useful comments and references. 
The work was partially supported by CONACyT grant
SEP-25750, and PROMEP grant UASLP-CA-21.}

From the works of Thompson \cite{thompson} we know  how the Jordan normal
form can be changed by a rank $k$ perturbation, see
Theorem~\ref{th3}. Particulary, it follows  that one can do everything
with a
geometrically simple spectrum by a rank $1$
perturbation, see Corollary~\ref{th1}. But the situation is quite
different if one restricts oneself to normal matrices, see
Theorem~\ref{th4} and Corollary~\ref{corol_Norm1}. We think that
Corollary~\ref{corol_Norm1} may be considered as a finite dimension
analogue of the continuous spectrum conservation under compact
perturbations in Hilbert spaces. For unitary and self-adjoint
matrices the inequality of Corollary~\ref{corol_Norm1} is the only
restrictions on "what can be done with a spectrum by a rank $k$
perturbation", see Theorem~\ref{Th_un_ad}. We don't know if there is an analogue
of Theorem~\ref{Th_un_ad} for normal matrices. It is worth to mention
that Corollary~\ref{corol_Norm1} for self-adjoint matrices follows from
Cauchy interlacing theorem \cite{Cauchy}. Theorem~\ref{Th_un_ad}
is related with the converse Cauchy interlacing theorem \cite{Fan}.

The spectrum of $H_1+H_2$ with known
spectra of self-adjoint matrices $H_1$ and $H_2$ is studied a lot, see  \cite{Klyachko}
and the bibliography therein.
Although the complete set of restrictions on the spectrum $H_1+H_2$
known in this situation, we are not sure that there is an easy proof
of Theorem~\ref{Th_un_ad} using results of \cite{Klyachko}.

Although Theorem~\ref{th3} should be known (see, for example, \cite{SSV2}, where
Theorem~\ref{th3} formulated in one direction), we will give a proof here, manly
because our proof falls in a general framework , which is also used in the proof
of Theorem~\ref{Th_un_ad}. Let us describe the framework. The set $\C_{n\times n}$
 of all complex $n\times n$-matrices (set of self-adjoint matrices) we equip with
 the arithmetic distance, $d(A,B)=\rank(A-B)$ (see \cite{Che}). The arithmetic distance is
geodesic for these cases. The spectral properties of matrices, such as Weyr characteristics
and spectra (multiset) also may be considered as a metric spaces with distance, related to
the arithmetic distance on matrices, see Section~\ref{sec_dist}. These distances also
turn out to be geodesic. Then we prove Theorem~\ref{th3} (Theorem~\ref{Th_un_ad})
for $\rank(A-B)=1$ and the general results will follow from Proposition~\ref{prop_geod}.
\begin{proposition}\label{prop_geod}
Let $X$ and $Y$ be geodesic metric spaces, let $O^X_n(x)$ denote the closed
ball of radius $n$ around $x$ in $X$,
let $\phi:X\to Y$ be such that $\phi(O^X_1(x))=O^Y_1(\phi(x))$ for all $x\in X$.
Then $\phi(O^X_n(x))=O^Y_n(\phi(x))$ for any $n\in \N$ and $x\in X$.
\end{proposition}
\begin{proof}
The proof is by induction. For $n=1$ there is nothing to prove. Step $n\to n+1$:
It follows that  $O^X_{n+1}(x)=\bigcup\limits_{z\in O^X_{n}(x)}O^X_1(z)$ ($X$ is geodesic),
then
$$
\phi(O^X_{n+1}(x))=\bigcup\limits_{z\in O^X_{n}(x)}\phi(O^X_1(z))=
\bigcup\limits_{z\in \phi(O^X_{n}(x))}O^Y_1(z)=\bigcup\limits_{z\in O^Y_{n}(\phi(x))}O^Y_1(z)=
O^Y_{n+1}(\phi(x)).
$$
\end{proof}

 In Section~\ref{sec_almost} and Section~\ref{sec_com} we use the
normalized arithmetic distance
$d_r(A,B)=\frac{\rank(A-B)}{n}$, where $n$ is the size of the
matrices. We are interested in the following questions: " Suppose
that matrices almost satisfy some equations (in the sense of
$d_r(\cdot,\cdot)$). If close to that matrices there exist matrices
satisfying the equations (uniformly with respect to  $n$)?" We
manage to answer only the following: close to an almost unitary
(self-adjoint) matrix there exists an unitary (self-adjoint) matrix.
We do not know if the same is true for normal matrices. (This
question has the affirmative answer for norm distance
$dn(A,B)=\|A-B\|$, see \cite{Lin}. It is equivalent to the
following: "close to any pair of almost commuting self-adjoint
matrices there exists a pair of commuting self-adjoint matrices (with
respect to the distance $dn(\cdot,\cdot)$). It is interesting that
there are almost commuting (with respect to $dn(\cdot,\cdot)$)
matrices, close to which there are no commuting matrices,
\cite{Choi,RuyTerry,Dan}). The similar question have been studied
for operators in Hilbert spaces (Calkin algebras, \cite{Farah}). In
Hilbert spaces the operator $a$ is called to be essentially normal
iff $aa^*-a^*a$ is a compact operator. In contrast with
Theorem~\ref{th6}, there exists an essentially unitary operator
which is not a compact perturbation of an unitary operator (just
infinite 0-Jordan cell). There is a complete characterization of
compact perturbations of normal operators, see \cite{Farah} and the
bibliography therein. Let us return to almost commuting matrices
with respect normalized arithmetic distance $d_r$.
In Section~\ref{sec_com}
we show that for any $A\in \C_{n\times n}$ with simple spectrum
there exists an almost commuting matrix, which is far from each
commuting with $A$ matrix. The similar problem for the pairs of
almost commuting matrices, as far as we know, is open.
Precisely, if for any $\epsilon>0$ there exists $\delta>0$ such that
for any $A,B$, $d_r(AB,BA)<\epsilon$ there exists $(\tilde A,\tilde B)$
satisfying $\tilde A\tilde B=\tilde B\tilde A$ and
$d_r(\tilde A,A),d_r(\tilde B,B)<\delta$
($\delta$ does not depends on size of the matrices).

We think that low rank perturbations of matrices may related to
sofic groups. The following question seems to be interesting from
this point of view (although it seems to goes beyond
the scope of the present article).
One can show that all solutions of equation
$C^{-1}A^{-1}CAC^{-1}AC=A^2$ in finite unitary matrices are trivial
in $A$ ($A=E$). On the other hand,  it is true that for any $\epsilon
>0$ there exist $A,C$, $d(A,E)=d(C,E)=1$ and
$d(C^{-1}A^{-1}CAC^{-1}AC,A^2)<\epsilon$. If the above assertion is
true with additional requirements $C^4=1$?
If not, it gives an example of non-sofic group.

Note. All linear spaces are supposed to be finite dimensional in the
rest of the article. $\C_{n\times n}$ will denote the set of all complex
$n\times n$-matrices, $\N=\{0,1,2,...\}$.

\section{Some discrete geodesic spaces.} \label{sec_dist}

\subsection{Arithmetic distance on $C_{n\times n}$}

\begin{lemma}\label{lm_arith}
The arithmetics distance $\rank(A-B)$ is geodesic on
\begin{itemize}
\item Set of all $n\times n$ matrices.
\item Set of all self-adjoint $n\times n$ matrices.
\item Set of all unitary $n\times n$ matrices.
\end{itemize}
\end{lemma}
\begin{proof}
It is clear that a rank $k$ matrix (self-adjoint matrix) may be represented
as sum of $k$ matrices (self-adjoint matrices) of rank $1$. The first two items
follow from the fact that set of matrices (self-adjoint matrices) is closed with
respect to summation.
For unitary matrices. Let $\rank(U_1-U_2)=k$, or, the same, $\rank(E-U_1^{-1}U_2)=k$.
It means that, in a proper basis, $U_1^{-1}U_2=\diag(\lambda_1,\lambda_2,...,\lambda_k,1,1,...,1)$.
Now the sequence $U_1,\;U_1\cdot\diag(\lambda_1,1,1,...,1),\;U_1\cdot
\diag(\lambda_1,\lambda_2,1,1,...,1)...
U_1\cdot\diag(\lambda_1,\lambda_2,...,\lambda_k,1,1,...,1)=U_2$
give us the geodesic needed.
\end{proof}
\begin{remark}
The methods used in the above proof are not applied for normal matrices -- the set
of normal matrices is not closed with respect neither summation nor multiplication.
In fact, an example from \cite{Fan}  hints that arithmetic distance might be non geodesic
on the set of normal matrices.
\end{remark}
\begin{proposition}\label{prop_mob1}
Let $\phi(x)=(ax+b)^{-1}(cx+d)$ be a M\"obius transformation of
$\C_{n\times n}$, defined on $A$, $B$. Then
$\rank(A-B)=\rank(\phi(A)-\phi(B))$.
\end{proposition}
\begin{proof}
A M\"obius transformation is a composition of linear transformations
$A\to aA+b$ ($a,b\in \C$)
and taking inverse $A\to A^{-1}$. Those transformations
(if defined) clearly conserve arithmetic distance,
for example, $\rank(A^{-1}-B^{-1})=\rank(A^{-1}(B-A)B^{-1})=\rank(A-B)$.
($A^{-1}$ and $B^{-1}$ is of full rank.)
\end{proof}

\subsection{Distance on  the spaces of the Weyr characteristics.}

Having in mind the Weyr characteristics of complex matrices (see below), we
introduce
the spaces $\Im_n$ of the Weyr characteristics.  Where $\Im_n$  is the space of
functions $\Z^+\times\C\to \N,\; (i,\la)\to\eta_i(\la)$ such that
\begin{itemize}
\item $\eta_i(\la)\neq 0$ for finitely  many $(i,\la)$ only, and
$\sum\limits_{\la\in\C}\sum\limits_{i\in\N}\eta_i(\la)=n$.
\item $\eta_i(\la)\geq \eta_{i+1}(\la)$.
\end{itemize}
On $\Im_n$  define a  metric
$d(\eta,\mu)=\max\limits_{(i,\la)}\{|\eta_i(\la)-\mu_i(\la)|\}$.
First of all let us note that $d(\cdot,\cdot)$ is indeed a metric.
Trivially, $d(\eta,\mu)=0$ implies $\eta=\mu$ and $d(\cdot,\cdot)$
satisfies triangle inequality as supremum (maximum) of semimetrics.
It is clear, that $d(\mu,\nu)$ is also well defined for $\mu$ and
$\nu$ in different spaces of Weyr characteristics (for different
$n$). We will need the following
\begin{proposition}\label{prop_weyr}
Let $\mu\in\Im_m$ and $n>m$. Then there exists $\nu\in\Im_n$ such
that for any $\eta\in\Im_n$,  the inequality $d(\mu,\eta)\geq
d(\nu,\eta)$ holds.
\end{proposition}
\begin{proof}
We can do as follows. Let $\mu_i(\lambda_0)\neq 0$ and
$\mu_{i+1}(\lambda_0)=0$. We can take
$\nu_{i+1}(\lambda_0)=\nu_{i+2}(\lambda_0)=...=\nu_{i+n-m}(\lambda_0)=1$
and $\mu_j(\lambda)=\nu_j(\lambda)$ for all other pairs
$(j,\lambda)$.
\end{proof}
\begin{proposition}
$\Im_n$ are  geodesic metric spaces.
\end{proposition}
\begin{proof}
Let $\eta,\mu\in \Im_n$ and $d(\eta,\mu)=k>1$ it is enough to find
$\nu\in\Im_n$ such that either $d(\eta,\nu)=1$ and $d(\nu,\mu)=k-1$,
or $d(\eta,\nu)=k-1$ and $d(\nu,\mu)=1$, moreover, by
Proposition~\ref{prop_weyr} it is enough to find $\nu\in\Im_m$ for
$m\leq n$. Let $S_+=\{(j,\lambda)\in
\Z^+\times\C\;|\;\eta_j(\lambda)-\mu_j(\lambda)=k\}$ and
$S_-=\{(j,\lambda)\in
\Z^+\times\C\;|\;\eta_j(\lambda)-\mu_j(\lambda)=-k\}$. Suppose, that
$|S_+|\geq |S_-|$ (if not, we can change $\eta\leftrightarrow\mu$).
Now let
$$
\nu_i(\lambda)=\left\{\begin{array}{lll}
           \eta_i(\lambda)-1 & \mbox{if} & (i,\lambda)\in S_+\\
           \eta_i(\lambda)+1 & \mbox{if} & (i,\lambda)\in S_-\\
           \eta_i(\lambda) & \mbox{if} & (i,\lambda)\not\in S_+\cup S_-
           \end{array}\right.
$$
We have to show that $\nu\in\Im_m$ for $m=n-|S_+|+|S_-|$. It is enough to
show that $\nu_{j+1}(\lambda)\leq \nu_j(\lambda)$. Suppose contrary, that
$\nu_{j+1}(\lambda)> \nu_j(\lambda)$. There are three possibility:
\begin{description}
\item[a)] $\eta_{j+1}(\lambda)=\eta_j(\lambda)$, $(j,\lambda)\in S_+$ and
$(j+1,\lambda)\not\in S_+$, but then
$k>\eta_{j+1}(\lambda)-\mu_{j+1}(\lambda)\geq \eta_{j+1}(\lambda)-\mu_j(\lambda)=
\eta_j(\lambda)-\mu_j(\lambda)=k$, a contradiction.
\item[b)] $\eta_{j+1}(\lambda)=\eta_j(\lambda)$, $(j+1,\lambda)\in S_-$ and
$(j,\lambda)\not\in S_-$, but then $-k<\eta_{j}(\lambda)-\mu_{j}(\lambda)\leq
\eta_{j}(\lambda)-\mu_{j+1}(\lambda)=
\eta_{j+1}(\lambda)-\mu_{j+1}(\lambda)=-k$, a contradiction.
\item[c)] $\eta_{j+1}(\lambda)=\eta_j(\lambda)-1$, $(j,\lambda)\in S_+$ and
$(j+1,\lambda)\in S_-$, but then
$-k=\eta_{j+1}(\lambda)-\mu_{j+1}(\lambda)\geq \eta_{j+1}(\lambda)-\mu_j(\lambda)=
\eta_j(\lambda)-\mu_j(\lambda)-1=k-1$, so $-k\geq k-1$ and $1/2\geq k$, a contradiction
with $k>1$.
\end{description}
Now, by construction, $d(\nu,\eta)=1$ and $d(\nu,\mu)=k-1$.
\end{proof}

\subsection{Distances $\dc$ and $\tilde\dc$ on finite multisets of the complex numbers.}

Multisets and operations. The language of multisets is very
convenient to deal with spectrums. We will need only finite
multisets. For a multiset $\A$ let $\set(\A)$ denote the set of
elements of $\A$, forgetting multiplicity. It is clear that
a muiltiset maybe considered as the  multiplicity function
$\chi_\A:\set(A)\to\N$, for any $x\not\in\A$ we will suppose
$\chi_\A(x)=0$. (For all cases, considered here,
$\set(\A)\subset\C$, so we can consider
$\chi_\A:\C\to\N=\{0,1,2...\}$.)  As far as the authors aware, there
are  several generalizations of the set-theoretical operations to
multisets. We will need the following operations:
\begin{itemize}
\item Difference of two multiset $\A\setminus\B$,
$\chi_{\A\setminus\B}(x)=\max\{0,\chi(\A)-\chi(\B))\}$.
\item Intersection $\A\cap X$ of a set $X$ and a multiset $\A$,
$$
\chi_{\A\cap X}(x)=\left\{\begin{array}{lll}
   \chi_\A(x),&\mbox{if} & x\in X \\
   0,&\mbox{if} & x\not\in X
   \end{array}\right.
$$
\item Union $\A\uplus\B$, $\chi_{\A\uplus\B}(x)=\chi_\A(x)+\chi_\B(x)$.
\end{itemize}
Let $S(a,r)=\{x\in\C\;: |x-a|\leq r\}$ and
$\S=\{S(a,r)\;|\;a\in\C,\;r\in \R^+\}$ be the set of circles. For
$\A,\B\subset_M\C$ let $\dc(\A,\B)=\max\limits_{S\in\S}\{||\A\cap
S|-|\B\cap S||\}$. Let as extend $\S$ to $\tilde \S$ which include
interior of complements of circles and semiplains: $\tilde S=S\cup
\{\{x\in\C\;:\;|x-a|\geq r\}\;|\;a\in\C,\;r\in\R^+\}\cup \{
\{x\in\C\;:\;\Imp(\frac{x-b}{a})\geq 0\}\;|\;a,b\in\C\}$; and introduce
new metric $\tilde\dc$:
$\tilde\dc(\A,\B)=\max\limits_{S\in\tilde\S}\{||\A\cap S|-|\B\cap
S||\}$.
\begin{proposition}
\begin{itemize}
\item $\dc$ and $\tilde\dc$ are metrics on the finite multisubsets of $\C$.
\item $\dc(\A,\B)=\dc(\A\setminus\B,\B\setminus\A)$,
$\tilde\dc(\A,\B)=\tilde\dc(\A\setminus\B,\B\setminus\A)$.
\item  If $|\A|=|\B|$, then $\tilde\dc(\A,\B)=\dc(\A,\B)$.
\end{itemize}
\end{proposition}
\begin{proof}
\begin{itemize}
\item The same as for the spaces of Weyr characteristics.
\item Let $\sum_S(\A,\B)=|\A\cap S|-|\B\cap S|=
\sum\limits_{x\in S}(\chi_\A(x)-\chi_\B(x))$.  Then
$\sum_S(\A\setminus\B,\B\setminus\A)=
\sum\limits_{x\in S}(\max\{0,\chi_\A(x)-\chi_\B(x)\}-\max\{0,\chi_\B(x)-\chi_\A(x)\})=
\sum\limits_{x\in S}(\chi_\A(x)-\chi_\B(x))$. Now the item follows by definition of
$\dc$ ($\tilde\dc$).
\item First of all, due to $\A$ and $\B$ are finite multisets, for any semiplain $p$
we can find a circle $c$ such that $\sum_p(\A,\B)= \sum_c(\A,\B)$.
Also for any closed circle $c_c$ there exists an open circle $c_o$ such that
$\sum_{c_c}(\A,\B)= \sum_{c_o}(\A,\B)$.
  Now, under
assumption of the item $\sum_S(\A,\B)=-\sum_{\C\setminus S}(\A,\B)$ and
the result follows.
\end{itemize}
\end{proof}

\begin{proposition}\label{prop_mob}
Let $\phi(x)=\frac{ax+b}{cx+d}$ be a M\"obius transformation of
$\C$, defined on $\set(\A)\cup\set(\B)$. Then
$\tilde\dc(\A,\B)=\tilde\dc(\phi(\A),\phi(\B))$.
\end{proposition}
\begin{proof}
A M\"obius transformation defines a bijection of $\tilde\S$.
\end{proof}

We don't know if the metric $\dc$ is geodesic on the multisets with fixed cardinality,
but its restriction on any circle o line is:

\begin{proposition}\label{lm_geod1}
Let $l\subset\C$ be a circumference or a straight line.  Let $\A,
\B\subset_M l$, $|\A|=|\B|=n$ and $\tilde\dc(\A,\B)=k\geq 2$. Then
there exists ${\cal C}\subset_M l$, $|{\cal C}|=n$ such that
$\tilde\dc(\A,{\cal C})=1$ and $\tilde\dc({\cal C},\B)=k-1$.
\end{proposition}
\begin{proof}
By Proposition~\ref{prop_mob}, it is enough to proof it for the unit circle.
Let us start with the case when $\set(\A)\cap\set(\B)=\emptyset$. Let
$\Gamma=\set(\A)\cup\set(\B)\subset C^1$. Let $|\Gamma|=r$.
We will cyclically anticlockwise
order $\Gamma=\{\gamma_0,\gamma_1,...,\gamma_{r-1}\}$ by elements of $Z_r$.
To construct $\cal C$ we move each element of $\A$
to the next element in $\Gamma$, precisely,
$\set({\cal C})\subseteq \Gamma$ and
$$
\chi_{\cal C} (\gamma_i)=\max\{0,\chi_\A(\gamma_i)-1\}+\chi_{\set(\A)}(\gamma_{i-1}),
$$
the other words
$$
\chi_{\cal C} (\gamma_i)=\left\{\begin{array}{lll}
                    \chi_\A(\gamma_i)-1 &\mbox{if} &
                      \gamma_i\in\set(\A)\;\mbox{and}\;\gamma_{i-1}\not\in\set(\A) \\
                    1 &\mbox{if} & \gamma_i\not\in\set(\A)\;\mbox{and}\; \gamma_{i-1}\in\set(\A) \\
                    \chi_\A(\gamma_i) & \mbox{for} & \mbox{the other cases}
                   \end{array}\right.
$$
We will check that $\cal C$ satisfies our needs.
For $x,y\in \Gamma$ let $[x,y]$ denote the closed segment of $C^1$, starting from $x$
and going anticlockwise to $y$ (so $[x,y]\cup [y,x]=C^1$).
It is clear that for $X,Y\subset_M\Gamma$ one has
$\tilde\dc(X,Y)=\max\{||X\cap [\alpha,\beta]|-|Y\cap [\alpha,\beta]||\;:\;\alpha,\beta\in\set(X)\cup\set(Y)\}$.
Denote by $\sum_{[\alpha,\beta]}(X,Y)=|X\cap [\alpha,\beta]|-|Y\cap [\alpha,\beta]|$.

Now, $\tilde\dc(\A,{\cal C})=\tilde\dc(\A\setminus{\cal C},{\cal C}\setminus\A)=1$, for
$\A\setminus{\cal C}$ and ${\cal C}\setminus\A$ are interlacing sets on $\C^1$.
Suppose further, $\tilde\dc({\cal C},\B)=\tilde\dc({\cal C}\setminus\B,\B\setminus{\cal C})\geq k$,
then there exists $[\gamma_i,\gamma_j]$ such that either
\begin{enumerate}
\item $\sum_{[\gamma_i,\gamma_j]}({\cal C}\setminus\B,\B\setminus{\cal C})
=\tilde\dc({\cal C}\setminus\B,\B\setminus{\cal C})\geq k$\\
or
\item $\sum_{[\gamma_i,\gamma_j]}({\cal C}\setminus\B,\B\setminus{\cal C})\leq -k$.
\end{enumerate}
In the first case we may assume that $\gamma_i,\gamma_j\in{\cal C}\setminus\B$
and $\gamma_{i-1},\gamma_{j+1}\not\in{\cal C}\setminus\B$. Now, changing
interval if necessary ($i^n=i-1$ and (or) $j^n=j-1$) we may, keeping
$\sum_{[\gamma_i,\gamma_j]}({\cal C}\setminus\B,\B\setminus{\cal C})$, achieve that
$\gamma_i,\gamma_j\in\A$
and $\gamma_{i-1},\gamma_{j+1}\not\in\A$. Then
$\sum_{[\gamma_i,\gamma_j]}(\A,\B)=
\sum_{[\gamma_i,\gamma_j]}({\cal C}\setminus\B,\B\setminus{\cal C})+1\geq k+1$,
a contradiction.
The second case may be considered similarly.

If $\set(\A)\cap\set(\B)\neq\emptyset$ then we can find ${\cal C'}$ for $\A\setminus\B$ and $\B\setminus\A$ and
then take ${\cal C}={\cal C'}\uplus X$, where $X=\A\setminus(\A\setminus\B)=\B\setminus(\B\setminus\A)$.
\end{proof}

\section{On spectrum of low rank perturbations.}
\begin{theorem}[Thompson]\label{tt}
Let $n \times n$ matrix $A$ over a field $F$ have similarity
invariants $h_n(A) \mid h_{n-1}(A) \mid  \ldots \mid h_1(A)$. Then:
as column $n$ tuple $x$ and row $n$-tuple $y$ range over all vectors
entries in $F$, the similarity invariants assumed by the matrix
$$
B=A+xy
$$
are precisely the monic polynomials $h_n(B) \mid \ldots
\mid h_1(B)$ over $F$ for which $degree(h_1(B) \cdots h_n(B))=n$ and
 $$
\begin{array}{l}
h_n(B) \mid h_{n-1}(A) \mid  h_{n-2}(B) \mid h_{n-3}(A) \mid \ldots,\\
h_n(A) \mid h_{n-1}(B) \mid  h_{n-2}(A) \mid h_{n-3}(B) \mid \ldots.
\end{array}
$$
\end{theorem}

We are going to reformulate Theorem~\ref{tt} for the field $\C$ using Weyr
characteristic.

Let $\eta_m (A, \lambda)$ denote the number of $\lambda$-Jordan
blocks in $A$ of size greater or equal to $m$ ($m \in \N$).
$$
\eta_m(A,\la)=\mbox{dim}\hspace{0.2cm} Ker(\la E - A)^m-\mbox{dim}\hspace{0.2cm} Ker(\la E - A)^{m-1}
$$
 This sequence of numbers $\eta_1(A,\la), \ldots \eta_q(A,\la)$
 is called the Weyr characteristic for the eigenvalue $\la$ of matrix $A$, see \cite{weyr}.

 \begin{theorem} \label{th3}
 Let $A\in \C_{n\times n}$  with Weyr invariants $\eta_m (A, \la)$.
 Then as $R$ ranges over all
 $n\times n$ complex matrices of rank less o equal $k$, the Weyr invariants
 assumed by the matrix $B=A+R$ are precisely those, that
satisfy both of the following conditions:
\begin{itemize}
\item For any $\la\in\C$ and any $m\in\N$
$$
\mid \eta_m (A, \la)-\eta_m (B, \la) \mid \leq k.
$$
\item $\sum\limits_{\la\in\C}\sum\limits_{m\in\N}\eta_m(B,\la)=n$.
\end{itemize}
\end{theorem}
\begin{proof}

It is enough to prove the theorem for $k=1$. Indeed, assume that the theorem is
valid for $k=1$. Then we may consider the
Weyr characteristics as a map $\psi:\C_{n\times n}\to \Im_n$, which satisfies
Proposition~\ref{prop_geod}. So, Theorem~\ref{th3} follows, though it states that
$\psi(O_k(A))=O_k(\psi(A))$.

Let us prove Theorem~\ref{th3} for
$k=1$.  For a given eigenvalue $\la \in sp(A)$ the sequence of
numbers
$$
q_1(A, \la) \geq q_2(A, \la) \geq ...,
$$
corresponding to the sizes of the $\la$-Jordan blocks in the Jordan
normal form of $A$ are know as the Segre characteristics of $A$
relative to $\la$, \cite{weyr}.

The similarity invariant factors of $A \in \C_{n \times n}$  are sequence of monic polynomials
 in
$x$, $h_n(A) \mid h_{n-1}(A) \mid  h_{n-2}(A) \ldots \mid
h_{1}(A)$. It is known that $h_i(A)=\prod_{\lambda} (\lambda -x)^{q_i(A,
\la)}$, where $\la \in sp(A)$ and $q_i(A, \la)$ is a Segre
characteristic corresponding to $\la$.
So, if $\rank(A-B)=1$, by Thompson's theorem \ref{tt} one has
\begin{equation}\label{ineq1}
\begin{array}{l}
q_1(B, \la) \geq q_2(A, \la) \geq q_3(B, \la) \geq \ldots ,\\
q_1(A, \la) \geq q_2(B, \la) \geq q_3(A, \la) \geq \ldots .
\end{array}
\end{equation}

As Weyr characteristic is the conjugate partition of Segre characteristic,
we can use Ferrers diagram to compute $\eta_m(B, \la)$ (see, \cite{weyr}) as the number of points
of column $m$ in Ferrer diagram of the Segre characteristic of $B$ relatively
to $\la$ (by short the Ferrer diagram of $(B, \la)$). Precisely,
the Ferrer diagram for $q(B,\la)$ is the set $F^\la_B=\{(i,j)\in\Z^+\times\Z^+\;|\;j\leq q_i(B,\la)\}$,
see Figure~\ref{figferrer} (the numbering from top to bottom and left to right).  The Weyr characteristics
is related with Ferrers diagram by the formula
$\eta_j(B,\la)=|\{(x,y)\in F^\la_B\;|\;x=j\}|$.
From inequalities (\ref{ineq1})
we have that $q_i(B, \la) \geq q_{i+1}(A, \la)$.
This inequality is equivalent to the statement
$\forall i\neq 1\; (i,j)\in F^\la_A \to (i-1,j)\in F^\la_B$
 (Figure~\ref{figferrer}), which is equivalent to the fact that
$\eta_j(B, \la) \geq \eta_j(A, \la)-1$. In similar form, from inequalities
(\ref{ineq1}) we can observe that $q_i(A, \la) \geq q_{i+1}(B, \la)$ with
$1 \leq i \leq n-1$. Therefore, $\eta_j(A, \la) \geq \eta_j(B, \la)-1$, and
the theorem follows.

\begin{figure}
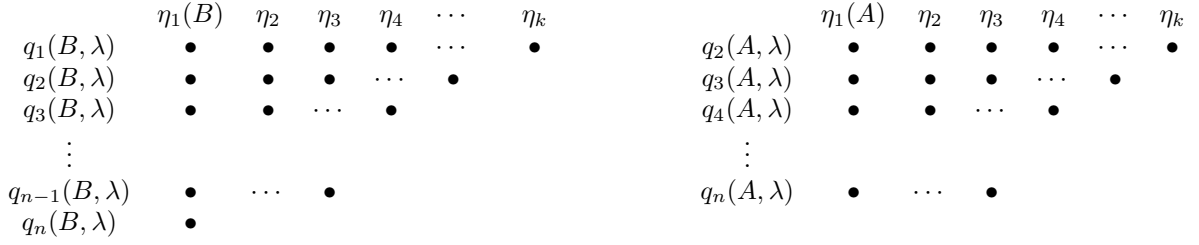

$$
\begin{array}{cccccccccccccccc}
\empty &    \eta_1(B) & \eta_2 &  \eta_3 &  \eta_4 & \cdots   &  \eta_k         &              &              &\eta_1(A) &   \eta_2&  \eta_3 & \eta_4  & \cdots &  \eta_k \\
q_1(B, \la) & \bullet & \bullet & \bullet & \bullet & \cdots  & \bullet         &\hspace{1.cm} & q_2(A, \la ) & \bullet  & \bullet & \bullet & \bullet & \cdots & \bullet\\
q_2(B, \la) & \bullet & \bullet & \bullet & \cdots & \bullet &
            \hspace{1cm} & & q_3(A, \la ) & \bullet  & \bullet & \bullet &\cdots & \bullet\\
 q_3(B ,\la) & \bullet  & \bullet & \cdots  & \bullet & &
            \hspace{1cm} & & q_4(A, \la ) & \bullet  & \bullet & \cdots & \bullet\\
 \vdots & & & & & & & &  \vdots\\
  q_{n-1}(B, \la) & \bullet  & \cdots  & \bullet  & & &
            \hspace{1cm} & & q_n(A, \la ) & \bullet  & \cdots & \bullet\\
    q_{n}(B, \la) & \bullet

\end{array}
$$
\caption{Relation between Ferrers diagrams of $(B, \la)$ and $(A,
\la)$} \label{figferrer}
\end{figure}
\end{proof}

\begin{corollP} \label{th1}
If the geometric multiplicity of any eigenvalue $\la$ of $A$ (number
of $\la$-Jordan cells)is $1$, then for any multiset $M$ of size $n$
there is a rank $1$ matrix $B$ such that $\spectr(A+B)=M$.
\end{corollP}

\section{Case of normal matrices.}

We will say that the vector $x$ is an $\alpha$-eigenvector if $Ax=\alpha
x$.  We will denote by $\RS(A,\lambda, \epsilon)$ the space
generated by all the $\alpha$-eigenvectors of $A$ with
$|\lambda-\alpha|\leq\epsilon$. For the case of normal matrices, the
following theorem shows that
 the difference between the dimension of  $\RS (A,\la,\e)$ and
 the dimension $\RS (B,\la,\e)$ is
 bounded by the rank of the difference matrix $A-B$.

\begin{theorem} \label{th4}
If $A$ and $B$ are normal matrices, then for any $\la$, and
for any $\e \geq 0$,
$$\mid dim( \RS (A,\la,\e)) - dim ( \RS (B,\la,\e)) \mid \leq \rank(A-B)$$
\end{theorem}

Let $X^\perp$ be the orthogonal complement of subspace $X$ and $P_X$
be an orthogonal projection on $X$.
\begin{lemma}\label{lm_ort}
Let $N:L\to L$ be a normal operator, and $X$ be a subspace of $L$
such that $\|(N-\la)x\|\leq\e\|x\|$ for any $x\in X$, then (we will
write $\RS (\la,\e)$ for $\RS (N,\la,\e)$)
\begin{enumerate}
\item $P_{\RS(\la,\e)}x\neq 0$ for any $x\in X$, $x\neq 0$.
\item $\|P_{\RS(\la,a\e)}x\|\geq \sqrt{1-\frac{1}{a^2}}\|x\|$ for
any $x\in X$.
\item $\dim(R(\la,\e))\geq\dim(X)$
\end{enumerate}
\end{lemma}
\begin{proof}
It is clear that (1) implies (3).

(1) Let $e_1,e_2,...,e_n$ be a diagonal orthonormal basis for $N$,
and $\la_1,\la_2,...\la_n$ corresponding eigenvalues
($Ne_i=\la_ie_i$). Let $x=\alpha_1e_1+\alpha_2e_2+...+\alpha_ne_n\in
X$ and $\|x\|=\sum\limits_{i=1}^n|\alpha_i|^2=1$. Let $x=x_1+x_2$
where $x_1\in \RS(\la,\e)$ and $x_2\in\RS^\perp(\la,\e)$. Now,
$\|(N-\la)x\|^2=\sum\limits_{i=1}^n|\alpha_i|^2|\la_i-\la|^2\leq\e^2$
implies that $\sum\limits_{i\;|\;|\la_i-\la|>\e}|\alpha_i|^2<1$. So,
$P_{\RS(\la,\e)}(x)=x_1=\sum\limits_{i\;|\;|\la_i-\la|\leq\e}\alpha_ie_i\neq
0$.

(2) Similarly,
$$
\sum\limits_{i\;|\;|\la_i-\la|>a\e}|\alpha_i|^2<\frac{1}{a^2}\;\;\mbox{and}\;\;
\sum\limits_{i\;|\;|\la_i-\la|\leq a\e}|\alpha_i|^2\geq
(1-\frac{1}{a^2}),
$$
so, $\|P_{\RS(\la,a\e)}x\|\geq \sqrt{1-\frac{1}{a^2}}\|x\|$
\end{proof}
Now we a ready to prove Theorem~\ref{th4}. Let $\rank(A-B)=r$
then there exists $X= \RS(A,\la,\e)\cap\ker(A-B)$ with
$\dim(X)\geq \dim(\RS(A,\la,\e))-r$ and $A|_X=B|_X$, so
$\|(B-\la)|_X\|\leq\e$ and, by Lemma~\ref{lm_ort},
$\dim(\RS(B,\la,\e))\geq \dim(\RS(A,\la,\e))-r$. By symmetry,
we get Theorem~\ref{th4}.

Theorem~\ref{th4} implies that any circle in complex plain,
containing $m$ spectral points of $A$ ($B$) should contain at least
$m-k$ spectral points of $B$ ($A$). So, we have
\begin{corollP}\label{corol_Norm1}
If $A$ and $B$ are  normal matrices then
$\dc(\spectr(A),\spectr(B))\leq \rank(A-B)$.
\end{corollP}
\begin{proof}
It is just a reformulation of Theorem~\ref{th4}.
\end{proof}

 If the condition of Corollary~\ref{corol_Norm1} describes all
accessible by rank $k$ perturbation spectra? We are going to show
that the answer is "yes" for self-adjoint and unitary matrices.
\begin{theorem}\label{Th_un_ad}
Let $A$ be a self-adjoint (unitary) $n\times n$-matrix. Let
$\B\subset_M\R$ ($\B\subset_M C^1$), $|\B|=n$. Then there exists
self-adjoint (unitary) matrix $B$ such that $\spectr(B)=\B$ and
$\rank(A-B)=\dc(\spectr(A),\B)$.
\end{theorem}
In fact the following, more general result is valid:
\begin{theorem}\label{Th_un_ad2}
Let $\l\subset\C$ be a circumference or straight line. Let $A$ be a
normal $n\times n$-matrix, $\spectr(A)\subset_M\l$. Let
$\B\subset_M\l$, $|\B|=n$. Then exists a normal matrix $B$ such that
$\spectr(B)=\B$ and $\rank(A-B)=\dc(\spectr(A),\B)$.
\end{theorem}
\begin{itemize}
\item It is enough to prove Theorem~\ref{Th_un_ad2}
for self adjoin matrices. Indeed, let $\spectr(A),\B\subset_M\l$, $|\B|=n$ for
a circle (line) $l\subset\C$. Then there exists a M\"obius transformation $\phi$,
defined on $\spectr(A)\cup\B$,
which map $l$ to the real line. Then $\phi(A)$ is a self-adjoint matrix and we
can apply Theorem~\ref{Th_un_ad} to $\phi(A)$ and $\phi(\B)$ to find
$\tilde B$ with $\spectr(\tilde B)=\phi(\B)$ and
$\rank(\phi(A)-\tilde B)=\tilde\dc(\phi(\spectr(A)),\phi(\B))$.
Now take $B=\phi^{-1}(\tilde B)$ and results follows, for the M\"obius
transformations conserve arithmetic distance on $C_{n\times n}$ and the distance
$\tilde\dc$ on multisets (Proposition~\ref{prop_mob1} and Proposition~\ref{prop_mob}).
\item It is enough to prove Theorem~\ref{Th_un_ad2}
for $\dc(\spectr(A),\B)=1$ and the rest will follow from Proposition~\ref{prop_geod},
Proposition~\ref{lm_geod1} and Lemma~\ref{lm_arith}.
\item Also w.l.g. we may assume that $\set(\spectr(A))\cap\set(\B)=\emptyset$. For
if $X=\spectr(A)\setminus(\spectr(A)\setminus\B)$ we can write
$A=A_1\oplus A_2$ with $\spectr(A_1)=X$ and
$\spectr(A_2)=\spectr(A)\setminus X$. We can find $B_2$ with
$\spectr(B_2)=\B\setminus X$ and $\rank(A_2-B_2)=1$.
Now, take $B=A_1\oplus B_2$.
\item Let $\A,\B\subset_M\R$, $\set(\A)\cap\set(\B)=\emptyset$, $|\A|=|\B|$ and
$\dc(\A,\B)=1$. Then, in fact, $\A$ and $\B$ are interlacing sets. It means that if
$\A=\{\alpha_1,\alpha_2,...,\alpha_n\}$,
$\B=\{\beta_1,\beta_2,...,\beta_n\}$ then $\alpha_1<\beta_1<\alpha_2<\beta_2<...$ or
$\beta_1<\alpha_1<\beta_2<\alpha_2<...$.
\end{itemize}

So, we need to prove only
\begin{lemma}\label{lm_rank1}
Let $A\in\C_{n\times n}$ be a self-adjoint matrix with a simple spectrum. Let
$\B\subset \R$ with $|\B|=n$. If $\spectr(A)$ and $\B$ are interlacing then
there exists a self-adjoint matrix $B$ with $\spectr(B)=\B$ and
$\rank(A-B)=1$.
\end{lemma}
Let $\spectr(A)=\A=\{\alpha_1,\alpha_2,...,\alpha_n\}$ and $\B=\{\beta_1,\beta_2,...,\beta_n\}$.
As $A$ and $B$ can be put in diagonal normal form
$\tilde A=\diag(\alpha_1,\alpha_2,...,\alpha_n)$ and $\tilde
B=\diag(\beta_1,\beta_2,...,\beta_n)$  by unitary transformations and unitary transformations
map (by conjugation) self-adjoint matrices to self adjoint matrices,
Lemma~\ref{lm_rank1} is equivalent to the fact that under our assumptions
on $\A$ and $\B$ the equation
\begin{equation} \label{eq_un}
\tilde A X- X\tilde B=R
\end{equation}
has a solution  in $(X,R)$ for unitary $X$ and $R$ of rank $1$.
Before solving Eq.\ref{eq_un} let us introduce some notations and prove a proposition.
For a finite $\A\subset\R$ let $P_\A(\lambda)=
\prod_{\alpha\in\A}(\lambda-\alpha)$. Let $\A$ and $\B$ be finite
subsets of $\R$ of equal cardinality. It follows from interpolation
that there exist unique $x:\A\to\A$, such that
\begin{equation}\label{eq_interpolation}
P_\B=P_\A-\sum_{\alpha\in\A}x_\alpha P_{\A\setminus\{\alpha\}}
\end{equation}
(we write $x_\alpha$ not $x(\alpha)$). Moreover,
$x_\alpha=P_\B(\alpha)/P_{\A\setminus\{\alpha\}}(\alpha)$.
Studying signs of $P_\B(\alpha)$ and $P_{\A\setminus\{\alpha\}}(\alpha)$ we
trivially get

\begin{proposition}\label{prop_interlace}
 If $\A$ and $\B$ interlacing, then all $x_\alpha$ in Eq.\ref{eq_interpolation}
have the same sign.
\end{proposition}
\begin{remark}
In fact inverse of this proposition is also valid, see
Lemma~1.20 of \cite{Fisk}.
\end{remark}
Now let us go back to  solutions of Eq.\ref{eq_un}. For $R=\{r_{ij}\}$ fixed the equation
has the unique solution $X=\{x_{ij}\}$, with
$x_{ij}=\frac{r_{ij}}{\alpha_i-\beta_j}$. Now,
suppose that $\rank(R)=1$ or, the same $r_{ij}=y_iz_j$ for some
$y,z\in\C^n$, $y,z\neq 0$. When the matrix $X$ is unitary? When its
columns (rows) are orthonormal, or
\begin{equation}\label{eq_unit}
z_jz^*_k\sum_{i}\frac{|y_i|^2}{(\alpha_i-\beta_j)(\alpha_i-\beta_k)}=\delta_{jk}
\end{equation}
It follows that $z_j\neq 0$ for all $j=1,...,n$, changing rows by
columns we get the same for $y$. So, the difficult part is to
guarantee that l.h.s. of Eq.~\ref{eq_unit} is $0$ for $j\neq k$.
Putting equality
$$
\frac{|y_i|^2}{(\alpha_i-\beta_j)(\alpha_i-\beta_k)}=\frac{1}{\beta_j-\beta_k}
(\frac{|y_i|^2}{\alpha_i-\beta_j}-\frac{|y_i|^2}{\alpha_i-\beta_k}).
$$
into Eq.\ref{eq_unit} and multiplying it by $P_\A(\beta_k)P_\A(\beta_j)(\beta_j-\beta_k)/(z_jz^*_k)$
 we get (after some elementary transformations):
$$
P_\A(\beta_k)\sum_{i=1}^n|y_i|^2P_{\A\setminus\{\alpha_i\}}(\beta_j)=
P_\A(\beta_j)\sum_{i=1}^n|y_i|^2P_{\A\setminus\{\alpha_i\}}(\beta_k),
$$
which imply that there exists $c$, such that
$$
\sum_{i=1}^n|y_i|^2P_{\A\setminus\{\alpha_i\}}(\beta_j)=cP_\A(\beta_j),
$$
for all $j$. Or the same, $\B$ is the set of roots of polynomial
$$
\Phi(\lambda)=\sum_{i=1}^n|y_i|^2P_{\A\setminus\{\alpha_i\}}(\beta_j)-cP_\A(\beta_j),
$$
so $|y_i|^2=P_\B(\alpha_i)/cP_{\A\setminus\{\alpha_i\}}(\alpha_i)$.
Choosing $c=1$ or $c=-1$, we get, by Proposition~\ref{prop_interlace}, that $|y_i|^2$ is well defined.
Now, take, for example, $y_i=|y_i|$. From Eq.\ref{eq_unit} for $j=k$ we can find $|z_j|^2$. Then,
taking $z_j=|z_j|$ we get needed solution of Eq.\ref{eq_un}.

\section{Almost unitary operators are near unitary operators with respect to normalized
arithmetic distance}\label{sec_almost}

For $A,B \in C_{n \times n}$, let $d_r(A,B)$ be the normalized arithmetic distance:
$$
d_r(A,B)=\frac{\rank(A-B)}{n}
$$

The matrix $A$ is called an $\alpha$-self-adjoint matrix if $d_r(A,
A^*)=\alpha$, where $A^*$ denotes the adjoint of $A$. The matrix $A$
is called an $\alpha$-unitary matrix if $d_r(A^*A,E) =\alpha$

The following theorems says that "near" to any $\alpha$-self-adjoint
matrix there exists a self-adjoint matrix $S$, and that "near" to
any $\alpha$-unitary matrix there exists an unitary matrix $U$ (for
small $\alpha$).
\begin{theorem} \label{th5}
For any $A \in \C_{n \times n}$ there exists a self-adjoint matrix $S$ ($S=S^*$) such that
$d_r(A,S) \leq d_r(A,A^*)$.
\end{theorem}
\begin{proof}
Take $S=\frac{1}{2}(A+A^*)$.
\end{proof}
\begin{theorem} \label{th6}
For any $A \in \C_{n \times n}$ there exists a unitary matrix $U$
($U^*U=E$), such that $d_r(A,U)\leq d_r(A^*A,E)$.
\end{theorem}
The good illustrations for this theorem are $0$-Jordan cells:
$$
  \left(
    \begin{array}{ccccccc}
      0 & 0 &  \cdots & 0 &0 & \cdots & 0\\
      1 & 0 &  \cdots  & 0 & 0& \cdots &0\\
         \multicolumn{7}{c}\dotfill\\
      0 & 0 &  \cdots & 0 & 0 & \cdots &0\\
      0 & 0 & \cdots &  1 & 0  & \cdots &0\\
         \multicolumn{7}{c}\dotfill\\

      0 & 0 & \cdots  & 0 & 0& \cdots & 0\\
    \end{array}
  \right)
  \left(
    \begin{array}{ccccccc}
      0 & 1 &  \cdots & 0 &0 & \cdots & 0\\
      0 & 0 &  \cdots  & 0 & 0& \cdots &0\\
         \multicolumn{7}{c}\dotfill\\
      0 & 0 &  \cdots &  0 & 1 & \cdots &0\\
      0 & 0 & \cdots &   0&  0  & \cdots &0\\
         \multicolumn{7}{c}\dotfill\\

      0 & 0 & \cdots  & 0 & 0& \cdots & 0\\
    \end{array}
  \right)=
\left(
    \begin{array}{ccccccc}
      0 & 0 &  \cdots & 0 &0 & \cdots & 0\\
      0 & 1 &  \cdots  & 0 & 0& \cdots &0\\
         \multicolumn{7}{c}\dotfill\\
      0 & 0 &  \cdots & 1 & 0 & \cdots &0\\
      0 & 0 & \cdots &  0 & 1 & \cdots &0\\
         \multicolumn{7}{c}\dotfill\\

      0 & 0 & \cdots  & 0 & 0& \cdots & 1\\
    \end{array}
  \right),
$$
but the matrix
$$
\left(
    \begin{array}{ccccccc}
      0 & 0 &  \cdots & 0 &0 & \cdots & 1\\
      1 & 0 &  \cdots  & 0 & 0& \cdots &0\\
         \multicolumn{7}{c}\dotfill\\
      0 & 0 &  \cdots & 0 & 0 & \cdots &0\\
      0 & 0 & \cdots &  1 & 0  & \cdots &0\\
         \multicolumn{7}{c}\dotfill\\

      0 & 0 & \cdots  & 0 & 0& \cdots & 0\\
    \end{array}
  \right)
$$
is unitary.
\begin{proof}
Let $\rank(A^*A-E)=r$, so there exist subspace $X\subset L$,
$\dim(X)=n-r$ such that $A^*A|_X=E|_X$. Consider $A|_X \,:\, X\to
Y=A(X)$. Under assumptions of the theorem $A^*(Y)=X$, it follows
that $(A|_X)^*=A^*|_Y:Y\to X$, so $A|_X\, :\,X\to Y$ is an unitary
operator. Choose any unitary operator $B: X^\perp\to Y^\perp$
($B^*B=E_{X^\perp}$). Then $U=A|_X\oplus B$ proves the theorem.
\end{proof}
It is not clear if this proof could be adapted for normal matrices
-- unitary operator from  an unitary space to another unitary space
is well defined, but how to define normal operators  between
different unitary spaces...?

Question: If we define $\alpha$-normal matrices in similar form to self-adjoint and unitary matrices,
the equivalent of theorems~\ref{th5} and~\ref{th6} are true for normal matrices?

\section{Almost commuting matrices}\label{sec_com}

\begin{theorem}\label{th_com}

For every $4 \leq n \in \N$ and every $A\in\C_{n \times n}$ with simple
spectrum there exists $X\in \C_{n\times n}$ such that
$d_r(AX,XA) < 2/n$  and for any matrix $B$, commuting with $A$, $d_r(B,X) \geq \frac{1}{2}$.

\end{theorem}

Before starting the proof of the theorem we need some facts.

\begin{proposition}\label{propol}
Let $\{\la_1,...,\la_k\}$ and $\{\al_1,...,\al_k\}$ be two disjoint sets,
then the matrix $M=[x_{ij}]$ with $x_{ij}=\frac{1}{\al_i-\la_j}$ is nonsingular
\end{proposition}
\begin{proof}
The matrix $M$ has the form
$$M=  \left(
          \begin{array}{cccc}
            \frac{1}{\al_1-\la_1} & \frac{1}{\al_1-\la_2} & \cdots & \frac{1}{\al_1-\la_k} \\
            \frac{1}{\al_2-\la_1}  & \frac{1}{\al_2-\la_2}  & \cdots & \frac{1}{\al_2-\la_k} \\            \multicolumn{4}{c}\dotfill\\
            \frac{1}{\al_k-\la_1} & \frac{1}{\al_k-\la_2} & \cdots & \frac{1}{\al_k-\la_k}    \\
                     \end{array}
        \right)
$$
Let $P(\la)=(\la-\la_1)(\la-\la_2)\cdots (\la-\la_k)$.  Multiply each row $j$ of the matrix $M$
by $P(\al_j)$ we obtain a matrix of the form
$$
\tilde M= \left(
          \begin{array}{cccc}
            P_1(\al_1) & P_2(\al_1) & \cdots & P_k(\al_1) \\
              P_1(\al_2) & P_2(\al_2) & \cdots & P_k(\al_2) \\
                 \multicolumn{4}{c}\dotfill\\
              P_1(\al_k) & P_2(\al_k) & \cdots & P_k(\al_k) \\
           \end{array}
        \right)
$$
where $P_j(\la)=\frac{P(\la)}{\la-\la_j}$.
This matrix will be nonsingular if and only if matrix $M$ is nonsingular.
We will prove that matrix $\tilde M$ is nonsingular showing that
the following system of linear  equations has a unique solution:
$$ \left(
          \begin{array}{cccc}
            P_1(\al_1) & P_2(\al_1) & \cdots & P_k(\al_1) \\
              P_1(\al_2) & P_2(\al_2) & \cdots & P_k(\al_2) \\
                 \multicolumn{4}{c}\dotfill\\
              P_1(\al_k) & P_2(\al_k) & \cdots & P_k(\al_k) \\
                     \end{array}
        \right)
        \left(
                 \begin{array}{c}
                   a_1\\
                   a_2 \\
                   \vdots \\
                   a_k \\
                 \end{array}
               \right)= \left(
                 \begin{array}{c}
                   b_1\\
                   b_2 \\
                   \vdots \\
                   b_k \\
                 \end{array}
               \right)
$$

So we have to solve the system $\sum_{i=1}^k a_iP_i(\al_j)=b_j$.
Consider the polynomial $\Phi(\la)=\sum_{i=1}^k a_iP_i(\la)$,
note that $\Phi(\la_i)=a_iP_i(\la_i)$ because $P_i(\la_j)=0$ for $i\neq j$.
We have $k$ points $b_j$, therefore, we can use Lagrange interpolation
to find the unique polynomial $\Phi(\la)$ of degree $k-1$ such that
$\Phi(\al_j)=b_j$, and then we can compute the values
$a_i=\frac{\Phi(\la_i)}{P_i(\la_i}$ ($P_i(\la_i)\neq 0$
for all $\la_i$ are different).
\end{proof}

Now we are ready to prove Theorem~\ref{th_com}.
\begin{proof}
Consider the matrix equation \begin{equation}\label{ax}
AX-XA=\{c_{ij}\},
\end{equation}
with $c_{ij}=i+j \textsf{mod} 2$. This matrix has the following form

$$AX-XA=  \left(
          \begin{array}{cccccc}
            0 & 1 & 0 & 1 & \cdots & 1+ n \hspace{0.2cm}\textsf{mod} 2 \\
            1 & 0 & 1 & 0 & \cdots &    n \hspace{0.2cm}\textsf{mod} 2 \\
            0 & 1 & 0 & 1 & \cdots & 1+n \hspace{0.2cm}\textsf{mod} 2 \\
            1 & 0 & 1 & 0  & \cdots & n \hspace{0.2cm}\textsf{mod} 2 \\
            \multicolumn{6}{c}\dotfill\\

               n +1  \hspace{0.2cm}\textsf{mod} 2 & \multicolumn{4}{c}\dotfill & 0
          \end{array}
        \right)
$$

Let $A$ be a diagonal matrix with simple spectrum $A=diag(\la_1,\la_2,...,\la_n)$.
One matrix $X$ that satisfies Eq.~\ref{ax} is $X = \{x_{ij}\}$  with \\
\begin{equation*}
x_{ij}=\begin{cases}

  \frac{c_{ij}}{\la_i-\la_j} & \text{for $i \neq j$} \\

  0 & \text{for $i=j$} \\
\end{cases}
\end{equation*}
Every matrix $B$ that commute with $A$ should be necessarily
a diagonal matrix $B=diag(b_1,b_2,...,b_n)$, then $X-B=\{x^*_{ij}\}$ with
\begin{equation*}
x^*_{ij}=\begin{cases}

  \frac{c_{ij}}{\la_i-\la_j} & \text{for $i \neq j$} \\

  -b_i & \text{for $i=j$} \\

\end{cases}
\end{equation*}

If we delete from $X-B$ the odd columns and the even rows, we obtain a
submatrix $X'$ of size $\lfloor \frac{n}{2} \rfloor \times \lfloor \frac{n}{2} \rfloor$ o
f the form
\begin{equation*}
\{x'_{ij}\}=  \frac{1}{\la_i^*-\la_j^*}
\end{equation*}
with $i^*=2i-1$, $j^*=2j$. By proposition~\ref{propol} this matrix is nonsingular.
Therefore, $\rank(X-B)\geq \frac{n}{2}$, and we obtain $d_r(X,B)\geq \frac{1}{2}.$
\end{proof}



\begin{thebibliography}{100}

 \bibitem{JLM}
 {\sc Atkinson, K.E.}
 \newblock An introduction to Numerical Analysus, 2nd ed.
 \newblock {\em John Wiley \& Sons, Ltd.}, 1989.

\bibitem{Cauchy}
{\sc Cauchy, A.L.}
\newblock Sur l'\'eqution \`a l'aide de laquelle on d\'etemine les in\'egalit\'es s\'eculaires,
\newblock {in Oeuvres Compl\`etes de A.L.Cauchy, 2nd Ser., V IX, Gauthier-Villars, Paris 1891, pp. 174-195}

\bibitem{Che}
{\sc Che-Hsien Wan, Zhe-Xian WAN, Lo-Keng Hua}
\newblock {\it Geometry of Matrices: In Memory of Professor L.K. Hua (1910-1985)}
\newblock World Scientific, 1996,
\newblock ISBN 9810226381, 9789810226381

  \bibitem{Choi}
    \newblock{\sc Choi, Man Duen,}
     \newblock Almost commuting matrices need not be nearly commuting,
   \newblock{\em Proc. Amer. Math. Soc.}, 102, (1988), N.3, pp. 529--533

 \bibitem{RuyTerry}
   \newblock{\sc Exel, R. and Loring, T.}
   \newblock Almost Commuting Unitary Matrices,
   \newblock {\it Proceedings of the American Mathematical Society}, Vol. 106, No. 4. (Aug., 1989), pp. 913-915.

\bibitem{Fan}
   \newblock {\sc Fan K. and Pall G.}
   \newblock Imbedding conditions for hermitian and normal matrix
   \newblock {\em Canada J. Math.}, 9 (1957) pp. 298-304

  \bibitem{Farah}
    \newblock {\sc Farah, Ilijas,}
    \newblock All automorphisms of the calkin algebra are inner,
     \newblock Preprint ArXiv:0705.3085v3

\bibitem{Fisk}
 \newblock {\sc Fisk, Steve,}
  \newblock Polynomials, roots, and interlacing,
  \newblock Preprint ArXiv: math/0612833

\bibitem{Klyachko}
\newblock {\sc Klyachko, Alexander,}
\newblock Vector Bundles, Linear Representations, and Spectral Problems,
\newblock Preprint ArXiv: math/0304325
 \bibitem{Lin}
 {\sc Lin, Huaxin},
 \newblock Almost commuting selfadjoint matrices and applications,
 \newblock in {\em Operator algebras and their applications} (Waterloo, ON,
              1994/1995), Fields Inst. Commun., 13, pp. 193--233, Amer. Math.
              Soc..


 \bibitem{JLFD}
 {\sc Moro, J. Dopico, F. M.}
 \newblock Low Rank Perturbation of Jordan Structure,
 \newblock {\em SIAM J. Matrix Anal. Appl.}, vol.~25, No 2, 2003, pp.~495--506.



 \bibitem{SSV}
 {\sc Savchenko, S. V.}
 \newblock On the change in the spectral properties of a matrix under a
              perturbation of a sufficiently low rank,
  \newblock {\em (Russian) Funktsional. Anal. i Prilozhen.}, 38 (2004), no. 1, pp.~85--88;
  translation in Funct. Anal. Appl. 38 (2004), no. 1, pp~69--71

\bibitem{SSV2}
 {\sc Savchenko, S. V.}
 \newblock Weyr characteristic under low rank perturbation,
  \newblock {\em Preprtin, L.D. Landau Institute for Theoretical Physics, Russian Academy of Science
Kosygina str. 2, Moscow, 119334, Russia}


\bibitem{weyr}
{\sc Shapiro, Helene,}
\newblock The Weyr characteristic.
\newblock {\em The American Mathematical Monthly Vol. 106, No. 10, \/} (1999), pp. 919--929.

\bibitem{thompson}
{\sc Thompson, Robertg}
\newblock Invariant factors under rank one perturbations.
\newblock {\it Can. J. Math.,}, XXXII, (1980), N.1. pp.240-245

  \bibitem{Dan}
   \newblock {\sc Voiculescu, Dan},
   \newblock Asymptotically commuting finite rank unitary operators without
              commuting approximants,
    \newblock {\it Acta Sci. Math. (Szeged)}, 45, (1983), N.1-4, pp. 429--431,


\end{thebibliography}
  \end{document}